\documentclass[11pt]{article}

\usepackage{amsmath}    
\usepackage{graphicx}   
\usepackage{verbatim}   
\usepackage{color}      
\usepackage{subfigure}  
\usepackage{amssymb,tikz}
\usepackage{soul}

\newtheorem{theorem}{Theorem}[section]
\newtheorem{proposition}[theorem]{Proposition}

\newtheorem{corollary}[theorem]{Corollary}
\newtheorem{lemma}[theorem]{Lemma}

\newcommand{\proof}{\noindent{\bf Proof.\ }}
\newcommand{\qed}{\hfill $\square$\medskip}
\def\cp{\,\square\,}
\def\sg{{\rm sg}}

\begin{document}

\title{Strong geodetic problem in grid like architectures}
 
\author{
	Sandi Klav\v zar $^{a,b,c}$
	\and
	Paul Manuel $^{d}$
}

\date{\today}

\maketitle
\begin{center}

	$^a$ Faculty of Mathematics and Physics, University of Ljubljana, Slovenia\\
	{\tt sandi.klavzar@fmf.uni-lj.si}\\
	\medskip
	
	$^b$ Faculty of Natural Sciences and Mathematics, University of Maribor, Slovenia\\
	\medskip
	
	$^c$ Institute of Mathematics, Physics and Mechanics, Ljubljana, Slovenia\\
	\medskip
	
	$^d$ Department of Information Science, College of Computing Science and Engineering, Kuwait University, Kuwait \\
	{\tt pauldmanuel@gmail.com}\\
	\medskip
		
\end{center}

\begin{abstract}
A recent variation of the classical geodetic problem, the strong geodetic problem, is defined as follows. If $G$ is a graph, then $\sg(G)$ is the cardinality of a smallest vertex subset $S$, such that one can assign a fixed geodesic to each pair $\{x,y\}\subseteq S$ so that these $\binom{|S|}{2}$ geodesics cover all the vertices of $G$. In this paper, the strong geodesic problem is studied on Cartesian product graphs. A general upper bound is proved on the Cartesian product of a path with an arbitrary graph and showed that the bound is tight on flat grids and flat cylinders.  
\end{abstract}


\medskip
\noindent{\bf AMS Subj.\ Class.: 05C12, 05C70, 05C76; 68Q17}

\section{Introduction}

Covering vertices of a graph with shortest paths is a problem that naturally appears in different applications; modelling them one arrives at different variations of the problem. The classical {\em geodetic problem}~\cite{HLT93} is to determine a smallest set of vertices $S$ of a given graph such that the (shortest-path) intervals between them cover all the vertices. The investigation on this problem up to 2011 has been surveyed in~\cite{BrKo11}, see also the book~\cite{pelayo-2013} for a general framework on the geodesic convexity. Recent developments on the geodetic problem include the papers~\cite{ekim-2014, ekim-2012, soloff-2015}, for a detailed literature survey see~\cite{MaKl16a,MaKl16b}. Another variation of the shortest-path covering problem is the {\em isometric path problem}~\cite{fisher-2001} where one is asked to determine the minimum number of geodesics required to cover the vertices, see also~\cite{pan-2006}. Motivated by applications in social networks, very recently the so called strong geodetic problem was introduced in~\cite{MaKl16a} as follows. 

Let $G=(V,E)$ be a graph. Given a set $S\subseteq V$, for each pair of vertices $\{x,y\}\subseteq S$, $x\ne y$, let $\widetilde{P}(x,y)$ be a {\em selected fixed} shortest path between $x$ and $y$. Then we set  
$$\widetilde{I}(S)=\{\widetilde{P}(x, y) : x, y\in S\}\,,$$ and let $V(\widetilde{I}(S))=\bigcup_{\widetilde{P} \in \widetilde{I}(S)} V(\widetilde{P})$. If $V(\widetilde{I}(S)) = V$ for some $\widetilde{I}(S)$, then the set $S$ is called a {\em strong geodetic set}. The {\em strong geodetic problem} is to find a minimum strong geodetic set $S$ of $G$. Clearly, the collection $\widetilde{I}(S)$ of geodesics consists of exactly $\binom{|S|}{2}$ paths. The cardinality of a minimum strong geodetic set is the {\em strong geodetic number} of $G$ and denoted by $\sg(G)$. For an edge version of the strong geodetic problem see~\cite{MaKl16b}.

In~\cite{MaKl16a} it was proved that the strong geodetic problem is NP-complete. Hence it is desirable to determine it on specific classes of graphs of wider interest. 
In this paper we follows this direction and proceed as follows. In the next section we first recall relevant properties of the Cartesian product of graphs. Afterwards we prove a lower bound on the strong geodetic number of Cartesian products in which one factor is a path. In Section~\ref{sec:dgrids-and-cylinders} we demonstrate that the bound is tight for grids $P_r\cp P_n$ and cylinders $P_r\cp C_n$ for the case when $r$ is large enough with respect to $n$; we will roughly refer to such graphs as {\em thin grids} and {\em thin cylinders}, respectively. But first we define concepts needed. 

All graphs considered in this paper are simple and connected. The {\em distance} $d_G(u,v)$ between vertices $u$ and $v$ of a graph $G$ is the number of edges on a shortest $u,v$-path alias $u,v$-{\em geodesic}. The {\em diameter} ${\rm diam}(G)$ of $G$ is the maximum distance between the vertices of $G$. We will use the notation $[n] = \{1,\ldots, n\}$ and the convention that $V(P_n) = [n]$ for any $n\ge 1$ as well as  $V(C_n) = [n]$ for any $n\ge 3$, where the edges of $P_n$ and of $C_n$ are defined in the natural way. 

The {\em Cartesian product} $G\cp H$ of graphs $G$ and $H$ is the graph with the vertex set $V(G) \times V(H)$, vertices $(g,h)$ and $(g',h')$ being adjacent if either $g=g'$ and $hh'\in E(H)$, or $h=h'$ and $gg'\in E(G)$. If $h\in V(H)$, then the subgraph of $G\cp H$ induced by the vertices of the form $(x,h)$,  $x\in V(G)$, is isomorphic to $G$;  it is denoted with $G^h$ and called a {\em $G$-layer}. Analogously $H$-layers are defined; if $g\in V(G)$, then the corresponding $H$-layer is denoted $H^g$. $G$-layers are also referred to as {\em horizontal layers} or, especially for grid as {\em rows}, while $H$-layers are {\em vertical layers} or {\em columns}.

\section{A lower bound on $\sg(P_r\cp G)$}
\label{sec:lower-bound}

In this section we prove a lower bound on the strong geodetic number of Cartesian products in which one factor is a path. For this sake we start by recalling some facts about the Cartesian product, especially about its metric properties. 

The Cartesian product is an associative and commutative operation. More precisely, the latter assertion means that the graphs $G\cp H$ and $H\cp G$ are isomorphic. We will implicitly (and explicitly) use this fact several times. Recall also that $G\cp H$ is connected if and only if both $G$ and $H$ are connected. Hence to assure that all graphs in this paper are connected, it suffices to assume that all the factor graphs are connected.   

The metric structure of Cartesian product graphs is well-understood, see~\cite[Chapter 12]{ikr-2008}. Its basis is the following result that was independently discovered several times, cf.~\cite[Lemma 12.1]{ikr-2008}. 

\begin{proposition}
\label{prop:distance-lemma}
If $(g,h)$ and $(g',h')$ are vertices of a Cartesian product $G \cp H$, then
$$d_{G \cp H}((g,h),(g',h')) = d_G(g,g') + d_H(h,h')\,.$$
\end{proposition}

If $(g,h)\in V(G\cp H)$, then the projections $p_G:V(G\cp H) \to V(G)$ and $p_H:V(G\cp H) \to V(H)$ are defined with $p_G((g,h)) = g$ and $p_H((g,h)) = h$. The projections $p_G$ and $p_H$ can be extended such that they also map the edges of $G\cp H$. More precisely, if $e=(g,h)(g',h)\in E(G\cp H)$, then $p_G(e) = gg'\in E(G)$, and if $e=(g,h)(g,h')\in E(G\cp H)$, then $p_G(e) = g\in V(G)$. Furthermore, we can also consider $p_G(X)$ and $p_H(X)$, where $X$ is a subgraph of $G\cp H$. 

Proposition~\ref{prop:distance-lemma} together with the fact that if $(g,h)$ and $(g',h)$ are vertices of the same $G$-layer, then every geodesic between $(g,h)$ and $(g',h)$ lies completely in the layer (see the first exercise in~\cite[12.3 Exercises]{ikr-2008}) implies the following. 

\begin{corollary}
\label{cor:projection}
Let $P$ be a geodesic in $G\cp H$. If $e=(g,h)(g',h)\in E(P)$, then $e$ is the unique edge of $P$ with $p_G(e) = gg'$. Moreover, $p_G(P)$ is a geodesic in $G$. 
\end{corollary}

Of course, by the commutativity of the Cartesian product, the assertions of Corollary~\ref{cor:projection} also hold for the projection of $P$ on $H$. 

After this preparation we can state the following technical lemma. 

\begin{lemma}
\label{lem:CP-extra-layer}
Let $G$ and $H$ be graphs, $\Omega$ be a minimum strong geodetic set of $G\cp H$, and $ \widetilde{I}(\Omega)$ its corresponding set of geodesics. If $|V(H)| > {\rm diam}(G)\binom{|\Omega|}{2} + |\Omega|$, then there exists a $G$-layer $G^{h}$ such that 
\begin{enumerate}
\item[(i)] $E(G^{h}) \cap \left( \cup_{P\in \widetilde{I}(\Omega)} E(P)\right) = \emptyset$ and 
\item[(ii)] $V(G^{h})\cap \Omega = \emptyset$.
\end{enumerate}
\end{lemma} 

\proof
Let  $t$ be the number of $G$-layers with the property that none of their edges lies on some path from $\widetilde{I}(\Omega)$, that is, 
$$t = \left| \left\{ h\in V(H):\ E(G^h) \cap \left( \cup_{P\in \widetilde{I}(\Omega)} E(P)\right) = \emptyset \right\}\right|\,.$$
Let $P$ be a geodesic from $\widetilde{I}(\Omega)$. By Corollary~\ref{cor:projection}, the edges of $P$ lie in at most ${\rm diam}(G)$ number of different $G$-layers. Hence, since $|\widetilde{I}(\Omega)| = \binom{|\Omega|}{2}$, the number of $G$-layers that contain edges of the paths from $\widetilde{I}(\Omega)$ is at most ${\rm diam}(G)\binom{|\Omega|}{2}$. Consequently, because we have assumed that $|V(H)| > {\rm diam}(G)\binom{|\Omega|}{2} + |\Omega|$ and as $|V(H)|$ is just the number of $G$-layers, we infer that $t > |\Omega|$. Therefore, by the pigeon-hole principle there exists at least one $G$-layer $G^{h}$, such that $E(G^{h}) \cap \left( \cup_{P\in \widetilde{I}(\Omega)} E(P)\right) = \emptyset$ and $V(G^{h})\cap \Omega = \emptyset$ as claimed. 
\qed

\begin{figure}[ht!]
\begin{center}
\scalebox{0.4}{\includegraphics{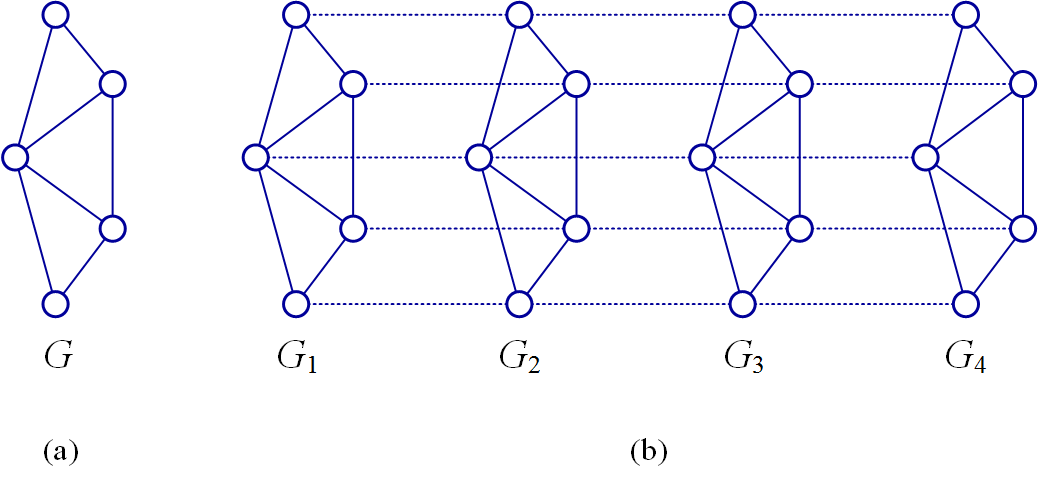}}
\end{center}
\caption{(a) Factor $G$. $(b)$ $P_4\cp G$, where the dotted edges are the edges of the $P_4$-layers, the other edges belong to the $G$-layers $G^1$, $G^2$, $G^3$, $G^4$.}
\label{fig:prod_graphs}
\end{figure}

We now restrict to Cartesian products of the form $P_r\cp G$. Since we have assumed that $V(P_n) = [n]$, the $G$-layers of $P_r\cp G$ are thus denoted with $G^1, \ldots, G^r$. See Fig.~\ref{fig:prod_graphs} for a graph $G$, the Cartesian product $P_4\cp G$, and the four $G$-layers. 

The main result of this section reads as follows. 

\begin{theorem}
\label{thm:LB_sgGxP3}
Let $\Omega$ be a minimum strong geodetic set of $P_r\cp G$. If $r > {\rm diam}(G)\binom{|\Omega|}{2} + |\Omega|$, then $\sg(P_r\cp G) \geq \lceil 2\sqrt{|V(G)|}\, \rceil$.
\end{theorem} 

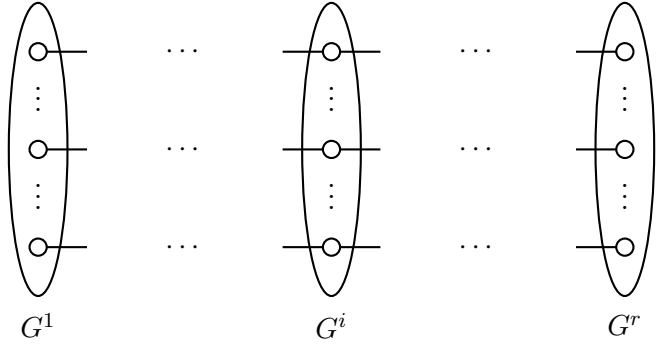
\begin{figure}[ht!]
\begin{center}
\begin{tikzpicture}[scale=0.65,style=thick]
\def\vr{5pt} 

\path (0,0) coordinate (a1);
\path (0,2) coordinate (a2);
\path (0,4) coordinate (a3);
\path (6,0) coordinate (a4);
\path (6,2) coordinate (a5);
\path (6,4) coordinate (a6);
\path (12,0) coordinate (a7);
\path (12,2) coordinate (a8);
\path (12,4) coordinate (a9);
\foreach \i in {1,2,...,6}
{ \draw (a\i) -- ++(1,0); }
\foreach \i in {4,5,...,9}
{ \draw (a\i) -- ++(-1,0); }

\foreach \i in {1,2,...,9}
{
   \draw(a\i)[fill=white] circle(\vr);
}
\foreach \i in {1,2,...,6}
{ \draw (a\i) ++(3,0)  node {$\cdots$}; }
\foreach \i in {1,2,4,5,7,8}
{ \draw (a\i) ++(0,1.2)  node {$\vdots$}; }
\draw (0,-1.6) node {$G^1$};
\draw (6,-1.6) node {$G^i$};
\draw (12,-1.6) node {$G^r$};
\draw (0,2) ellipse (0.6cm and 3cm);
\draw (6,2) ellipse (0.6cm and 3cm);
\draw (12,2) ellipse (0.6cm and 3cm);
\end{tikzpicture}
\end{center}
\caption{The $G$-layer $G^i$ guaranteed by Lemma~\ref{lem:CP-extra-layer} is neither the first nor the last $G$-layer.}
\label{fig:prod_graphs_partition}
\end{figure}

\proof
Applying Lemma~\ref{lem:CP-extra-layer} we infer that $P_r\cp G$ contains a $G$-layer, say $G^{i}$, such that no edge of $G^i$ lies on paths from $\widetilde{I}(\Omega)$ and such that $V(G^i)\cap \Omega = \emptyset$. Note that $i\ne 1$ and $i\ne r$, for otherwise the vertices of $G^1$ (resp.\ $G^r$) would not be covered with the paths from $\widetilde{I}(\Omega)$. We can hence partition $\Omega$ into non-empty sets $\Omega_1$ and $\Omega_2$ by setting
$$\Omega_1 = \Omega \cap \left(\bigcup_{j=1}^{i-1}G^j\right)\quad {\rm and}\quad
\Omega_2 = \Omega \cap  \left(\bigcup_{j=i+1}^{r}G^j\right)\,,$$
cf.~Fig.~\ref{fig:prod_graphs_partition}.

In order to cover the vertices of $G^i$, we must have $|\widetilde{I}(\Omega)| \geq |V(G)|$. This in turn implies that $|\Omega_1|\cdot |\Omega_2| \ge |V(G)|$. Since the arithmetic mean is not smaller than the geometric mean, we have 
$$\frac{|\Omega|}{2} = \frac{|\Omega|_1 + |\Omega_2|}{2} \ge \sqrt{|\Omega_1|\cdot |\Omega_2|} \ge \sqrt{|V(G)|}\,.$$
Thus $|\Omega|\geq 2\sqrt{|V(G)|}$. Since the number of vertices is an integer we conclude that $\sg(P_r\cp G) \geq \lceil 2\sqrt{|V(G)|}\, \rceil$.  
\qed

\section{Thin grids and cylinders}
\label{sec:dgrids-and-cylinders}

In this section we determine the strong geodetic number of thin grids $P_r\cp P_n$ ($r\gg n)$ and thin cylinders $P_r\cp C_n$ ($r\gg n)$. (The geodetic number in Cartesian products was investigated in~\cite{BrKl08, soloff-2015}, while in~\cite{caseras-2010} and in~\cite{bresar-2011} it was studied on strong products and lexicographic products, respectively.) Recall that by our convention on the vertex sets of paths and cycles, $V(P_r\cp P_n) = V(P_r\cp C_n) = \{(i,j):\ i\in [r], j\in [n]\}$.  

\begin{lemma}
\label{LUBWGS}
If $2 \le n \leq r$, then $\sg(P_r\cp P_n)\leq \lceil 2\sqrt{n}\, \rceil$.
\end{lemma} 

\proof
In order to prove the inequality, we need to construct a strong geodetic set of cardinality $\left\lceil 2 \sqrt{n}\, \right\rceil$. 

We first consider the case when $n$ is a perfect square, $n=k^2$. For each $i\in [k]$ define the vertices $a_i$ and $b_i$ of $P_r\cp P_n$ with  
\begin{quote}
$a_i = (1,(i-1)k+i)$, \\
$b_i = (r,(i-1)k+i)$,
\end{quote}
and set $S=\{a_1,a_2, \ldots, a_k\} \cup \{b_1,b_2, \ldots, b_k\}$. Note that $a_1=(1,1)$, $b_1=(r,1)$, $a_k=(1,k^2)$, and $b_k=(r,k^2)$ are the four vertices of $P_r\cp P_n$ of degree $2$ and that $|S| = 2k = 2\sqrt n$. Now we show that $S$ is a strong geodetic set of $P_r\cp P_n$ by constructing $\widetilde{I}(S)$ such that  $V(\widetilde{I}(S)) = V(P_r\cp P_n)$. It will suffice to select a geodesic for each pair of vertices $a_i$ and $b_j$ to achieve our goal.  

\begin{figure}[htb!]
	\begin{center}
	\scalebox{0.47}{\includegraphics{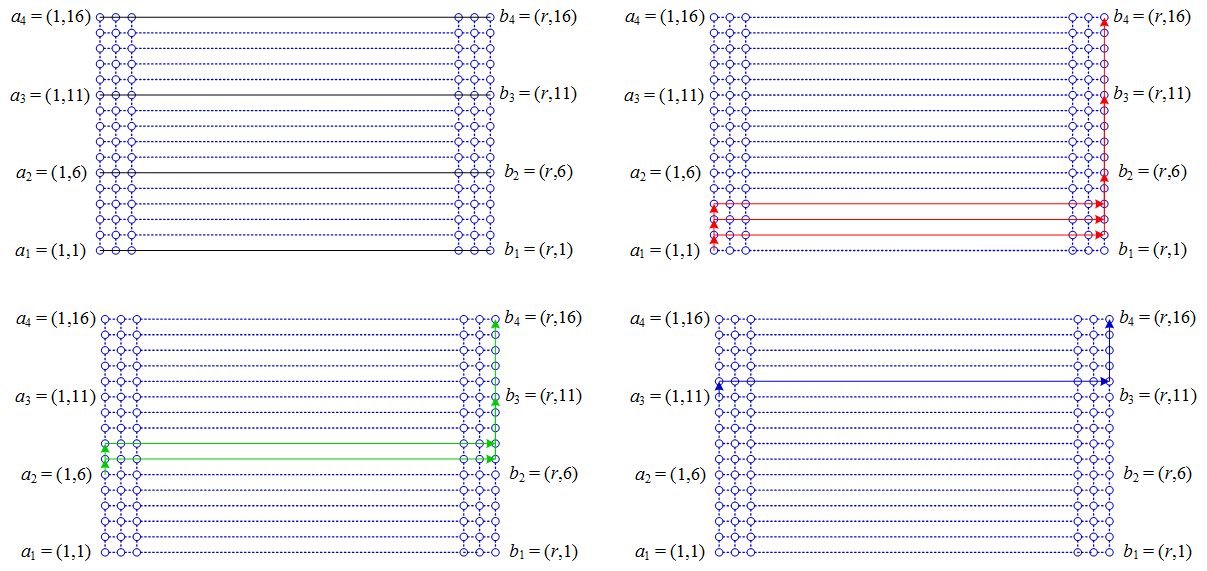}}
	\end{center}
	\caption{$P_r\cp P_{4^2}$ and the geodesics from $\widetilde{I}(S)$ between $a_i$ and $b_j$, where $i < j$.}
	\label{fig:sg_UB_grid_r2x16}
\end{figure}

For each $i\in [k]$ there is a unique $a_i,b_i$-geodesic which thus must belong to $\widetilde{I}(S)$. We next inductively add geodesics to $\widetilde{I}(S)$. First, add geodesics to $\widetilde{I}(S)$ one by one  from $a_1$ to respectively $b_2, \ldots, b_k$ as follows. Start in $a_1$ and traverse the first column until the first not yet traversed row is reached. Then traverse this row and complete the path by traversing the last column until the vertex $b_j$ that is just considered is reached. These paths thus cover $k-1$ rows. Then proceed along the same way for the vertices $a_2, \ldots, a_{k-1}$, respectively covering $k-2, \ldots, 1$ new rows. Repeat next the above procedure by inductively construction the geodesics from $a_k, \ldots, a_2$ to the vertices $b_j$. In this way the remaining rows are covered. The construction is illustrated in Fig.~\ref{fig:sg_UB_grid_r2x16}.

Assume next that $n = k^2 + \ell$, where $1\le \ell \le k$. In this case $\left\lceil 2\sqrt{n}\,\right\rceil = 2k+1$. Define the vertices $a_i$ and $b_i$ as above and set 
$$S = \{a_1,a_2, \ldots, a_k\} \cup \{b_1,b_2, \ldots, b_k\} \cup \{(1,n)\}\,.$$
As in the previous case, all the vertices of the subgraph $P_{r}\cp P_n$ can be covered using the vertices from $\{(1,i^2),(r,i^2):\ i\in [k]\}$. Then it is not difficult to show that the vertex $(1,n)$ will take care of the remaining vertices of $P_r\cp P_n$. Since $|S| = 2k + 1$, we are done also in this case. 

Finally, suppose that $n = k^2 + \ell$, where $k+1\le \ell \le 2k$. In this case $\left\lceil 2\sqrt{n}\,\right\rceil = 2k+2$. Setting 
$$S = \{a_1,a_2, \ldots, a_k\} \cup \{b_1,b_2, \ldots, b_k\} \cup \{(1,n),(r,n)\}$$
we can argue similarly as above that $S$ is a strong geodetic set of $P_r\cp P_n$.
\qed

The first main result of this section reads as follows. 

\begin{theorem}
\label{WGN_grids}
If $r > \binom{\left\lceil 2\sqrt{n}\right\rceil}{2}(n-1)+\left\lceil 2\sqrt{n}\right\rceil$, then $\sg(P_r\cp P_n)= \lceil 2\sqrt{n}\, \rceil$. 
\end{theorem}
\proof
Let $\Omega$ be a minimum strong geodetic set of $P_r\cp P_n$. By Lemma~\ref{LUBWGS} we know that $|\Omega| \leq \lceil 2\sqrt{n}\, \rceil$. Hence 
$$r > \binom{\left\lceil 2\sqrt{n}\right\rceil}{2}(n-1)+\left\lceil 2\sqrt{n}\right\rceil 
\ge \binom{|\Omega|}{2}(n-1)+|\Omega|\,.$$
As ${\rm diam}(P_n) = n-1$, Theorem~\ref{thm:LB_sgGxP3} implies that $\sg(P_r\cp P_n) = |\Omega| \geq \lceil 2\sqrt{n}\, \rceil$.
\qed

\begin{figure}[htb!]
	\begin{center}
		\scalebox{0.7}{\includegraphics{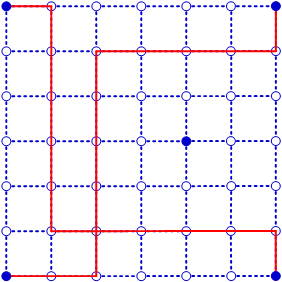}}
	\end{center}
	\caption{The blue bullets form a strong geodetic set of $P_7\cp P_7$. For the sake of clarity, only two geodesics are drawn. The reader can easily identify other geodesics.}
	\label{fig:Grids_Count_Ex_LB}
\end{figure}

Of course, it would be desirable to determine the exact strong geodetic number for all grids. To see that Theorem~\ref{WGN_grids} cannot be extended to all grids consider the product $P_7\cp P_7$. In Fig.~\ref{fig:Grids_Count_Ex_LB} we have produced a strong geodetic set consisting of 5 vertices. Thus $\sg(P_7\cp P_7)\leq 5 < \left\lceil 2\sqrt{7}\right\rceil = 6$. 

In the rest of the section we consider cylinders $P_r\cp C_n$ and prove an analogous result for thin cylinders ans we did for thin grids. We start with: 

\begin{lemma}
	\label{LLBWGS2_cyl}
	If $2 \le n \leq r$, then $\sg(P_r\cp C_n)\leq \lceil 2\sqrt{n}\, \rceil$.	
\end{lemma}

\proof
The proof of Lemma \ref{LUBWGS} is modified to accommodate cylinders. We first consider the case when $n$ is a perfect square, $n=k^2$. For each $i\in [k]$ define the vertices $a_i$ and $b_i$ of $P_r\cp C_n$ with  
\begin{quote}
$a_i = (1,(i-1)k+1)$, \\
$b_i = (r,(i-1)k+1)$,
\end{quote}
and set $S=\{a_1,a_2, \ldots, a_k\} \cup \{b_1,b_2, \ldots, b_k\}$. We claim that $S$ is a strong geodetic set of $P_r\cp C_n$ by constructing $\widetilde{I}(S)$ such that  $V(\widetilde{I}(S)) = V(P_r\cp C_n)$. 

\begin{figure}[htb!]
	\begin{center}
		\scalebox{0.43}{\includegraphics{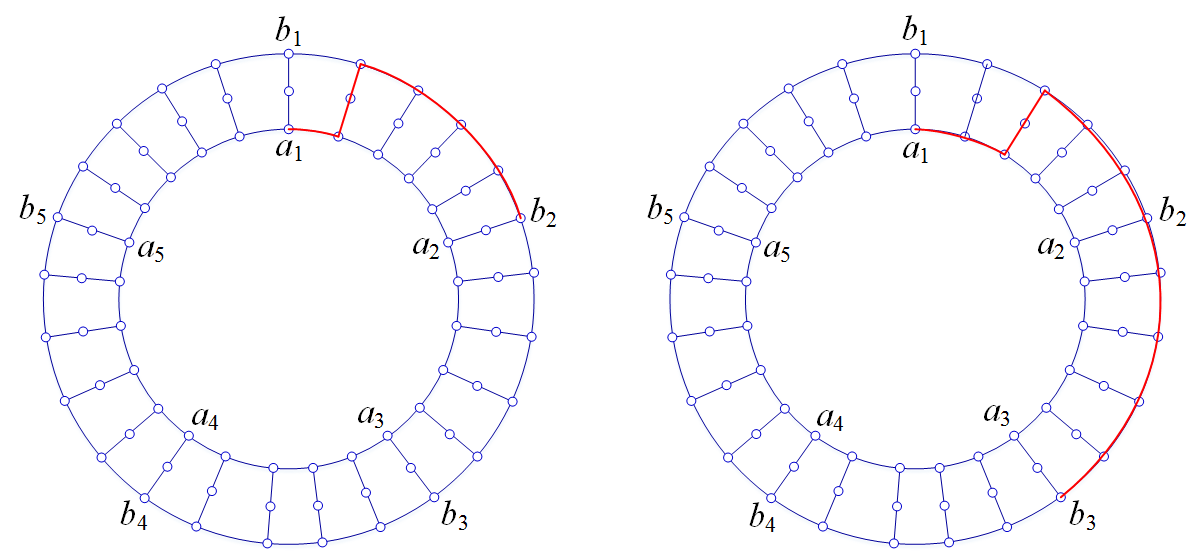}}
	\end{center}
	\caption{The paths marked by red edges are members of $\widetilde{I}(S)$. The first geodesic is between $a_1$ and $b_2$ and the second geodesic is between $a_1$ and $b_3$.}
	\label{fig:sg_UB_cylinder3x25_a}
\end{figure}

\begin{figure}[htb!]
	\begin{center}
		\scalebox{0.43}{\includegraphics{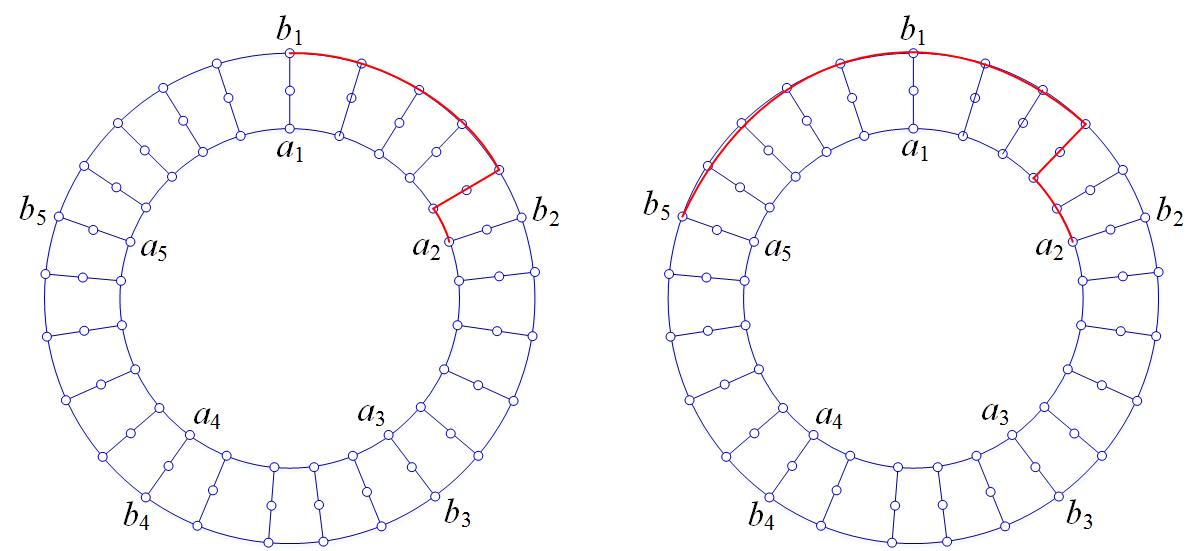}}
	\end{center}
	\caption{The paths marked by red edges are members of $\widetilde{I}(S)$. The first geodesic is between $a_2$ and $b_1$ and the second geodesic is between $a_2$ and $b_5$.}
	\label{fig:sg_UB_cylinder3x25_b}
\end{figure}

Assume first that $k$ is odd. Select the geodesics between $a_1$ and $b_2, b_3, \ldots, b_{\lfloor k/2\rfloor}$ and geodesics between $a_2$ and $b_1, b_{k}, \ldots, b_{\lceil k/2\rceil+2}$ such that all the $P_r$-layers $P_r^2, \ldots, P_r^{k-1}$ are covered by them. See Figs.~\ref{fig:sg_UB_cylinder3x25_a} and~\ref{fig:sg_UB_cylinder3x25_b}, where the case $P_2 \cp C_{25}$ is illustrated; that is, $r=2$ and $k=5$.     

By symmetry we can design geodesics starting from $a_i$ and $a_{i+1}$ such that all the corresponding $P_k$-layers are covered. In conclusion, $S$ is a strong geodetic set. If $k$ is even, we can proceed similarly to verify that $S$ is a strong geodetic set also in this case. Finally, if $n = k^2 + \ell$ and $1\le \ell \le 2k$, then by a similar construction as in the proof of Lemma~\ref{LUBWGS} (by adding either one or two new vertices to $S$, depending on whether $1\le \ell \le k$ or $k+1\le \ell \le 2k$, the construction is completed.  
\qed

We are now ready for the second main result of this section. It will be deduced similarly as Theorem~\ref{WGN_grids}.

\begin{theorem}
\label{WGN_cylinders}
If $r > \binom{\left\lceil 2\sqrt{n}\right\rceil}{2} \left\lfloor \frac{n}{2} \right\rfloor + \left\lceil 2\sqrt{n}\right\rceil$, then $\sg(P_r\cp C_n)= \lceil 2\sqrt{n}\, \rceil$. 
\end{theorem}

\proof
Let $\Omega$ be a minimum strong geodetic set of $P_r\cp C_n$. From Lemma~\ref{LLBWGS2_cyl} we know that $|\Omega| \leq \lceil 2\sqrt{n}\, \rceil$. Hence 
$$r > \binom{\left\lceil 2\sqrt{n}\right\rceil}{2} \left\lfloor \frac{n}{2} \right\rfloor + \left\lceil 2\sqrt{n}\right\rceil 
\ge \binom{|\Omega|}{2} \left\lfloor \frac{n}{2} \right\rfloor + |\Omega|\,.$$
As ${\rm diam}(C_n) = \left\lfloor \frac{n}{2} \right\rfloor$, Theorem~\ref{thm:LB_sgGxP3} implies that $\sg(P_r\cp C_n) = |\Omega| \geq \lceil 2\sqrt{n}\, \rceil$.
\qed

\section{Further research}
\label{sec:Fur-Research}

In this paper we have studied the strong geodetic problem on Cartesian product graphs and determined the strong geodetic number for ``flat" grids and cylinders. The first natural problem is of course to determine the strong geodetic number for all grids and cylinders. Next it would be interesting to consider the strong geodetic number on additional interesting Cartesian product graphs, such as torus graphs (product of two cycles) as well as on general Cartesian products. More generally, we can ask for the strong geodetic number on Cartesian product of more than two factors, in particular on multidimensional grid graphs.

\section*{Acknowledgement}

S.K.\ acknowledges the financial support from the Slovenian Research Agency
(research core funding No.\ P1-0297).

\end{document}